\documentclass{article}

\usepackage{enumerate}
\usepackage{amsmath,xspace,amssymb,mathrsfs}

\input xy
\xyoption{all}
\CompileMatrices

\title{Note on the construction of free monoids}
\author{Stephen Lack 
\thanks{The support of the Australian Research Council and
DETYA is gratefully acknowledged.}
\\School of Computing and Mathematics\\
University of Western Sydney\\
Locked Bag 1797 Penrith South DC NSW 1797\\
Australia\\
email: {\tt s.lack@uws.edu.au}}
\date{}

\newcommand{\C}{{\ensuremath{\mathscr C}}\xspace}
\newcommand{\E}{{\ensuremath{\mathscr E}}\xspace}
\newcommand{\M}{{\ensuremath{\mathscr M}}\xspace}

\newcommand{\PP}{{\ensuremath{\mathbb P}}\xspace}

\newcommand{\MonC}{\textnormal{Mon\C}\xspace}
\newcommand{\PtC}{\textnormal{Pt\C}\xspace}

\newcommand{\DD}{{\ensuremath{\mathbf\Delta}}\xspace}
\newcommand{\DDm}{{\ensuremath{\mathbf{\Delta_{mon}}}}\xspace}

\newcommand{\Set}{\textnormal{\bf Set}\xspace}
\newcommand{\Mon}{\textnormal{\bf Mon}\xspace}
\newcommand{\Grp}{\textnormal{\bf Grp}\xspace}
\newcommand{\Span}{\textnormal{\bf Span}\xspace}
\newcommand{\Prof}{\textnormal{\bf Prof}\xspace}
\newcommand{\Cat}{\textnormal{\bf Cat}\xspace}

\newcommand{\ot}{\otimes}
\newcommand{\colim}{\textnormal{colim}\xspace}

\newtheorem{theorem}{Theorem}
\newtheorem{lemma}[theorem]{Lemma}
\newtheorem{proposition}[theorem]{Proposition}

\begin{document}

\label{firstpage}
\maketitle

\begin{abstract}
We construct free monoids in a monoidal category $(\C,\ot,I)$ with
finite limits and countable colimits, in which tensoring on either side
preserves reflexive coequalizers and colimits of countable chains. In
particular this will be the case if tensoring preserves sifted colimits.
\end{abstract}

\section{Background}

For any monoidal category $(\C,\ot,I)$, one can form the category of 
monoids in \C, and for suitable choice of \C, this contains many important
notions, such as monoids, rings, categories, differential graded algebras,
and monads: see \cite[Chapter~VII]{CWM}. For each such \C, the category
\MonC of monoids in \C has a forgetful functor $U:\MonC\to\C$, and this
forgetful functor often has a left adjoint, sending an object of \C to the
free monoid on that object. In particular, if \C has countable coproducts,
and these are preserved by tensoring on either side, then the free monoid
on $X$ is given by the well-known ``geometric series''
$$I+X+X^2+X^3+\ldots$$
where $X^n$ stands for the $n$th``tensor power'' $X\ot\ldots\ot X$ of $X$. 

This case includes the
free monoid, the free ring (on an additive abelian group), the free category 
(on a graph) and the free differential graded algebra (on a chain complex),
but it does not help with the case of free monads. Conditions for the 
existence of free monads were given by Barr in \cite{Barr:free-triples}.
Further analysis of the free monoid construction was given by 
Dubuc in \cite{Dubuc:free-monoids}, of which more will be said below.

The epic paper \cite{Kelly-transfinite} of Kelly analyzes many constructions of 
free monoids, free algebras, and colimits, generally requiring transfinite
processes. It provides very general conditions for the existence of free
monoids, in the case where the category \C is cocomplete and the functors
$-\ot C:\C\to\C$ are cocontinuous, for all objects $C$ (the conditions on 
$C\ot-$ are then quite mild). This allows the construction of free monads
in many cases, and also the construction of free operads \cite{May-operads,
Kelly-operads}. For example if each $C\ot-$ preserves filtered colimits
then the free monoid is given by a ``factorized'' version of the geometric
series:
$$I+X(I+X(I+X(I+\ldots$$

A recent paper of Vallette \cite{Vallette:free-monoids} gave a construction
of free monoids under much stronger assumptions on $C\ot-$ but weaker 
assumptions on $-\ot C$, and moreover under the assumption that \C is abelian.
These assumptions allowed the construction of free properads 
and other closely related free structures, in the abelian context. 

In this paper
we generalize and simplify the construction of Vallette, removing the assumption
that \C is abelian. Specifically, we suppose that 
\C has finite limits and countable colimits, and that the functors $-\ot C$ and 
$C\ot -$ preserve reflexive coequalizers and colimits of countable chains, 
and we construct free monoids under these assumptions. (This would be the 
case for example if tensoring preserved sifted colimits \cite{sifted}.)
Some examples of such monoidal categories
are given in Section~\ref{sect:examples}. The question of whether these
free monoids are algebraically free \cite{Kelly-transfinite} is briefly
discussed in Section~\ref{sect:algebraically-free}.

\subsection*{Notation}

The tensor product of objects will generally be denoted by juxtaposition:
$XY$ stands for $X\ot Y$ (just as $X^2$ stands for $X\ot X$). We sometimes 
write as if the monoidal structure on \C were strict. This is merely for
convenience; by the coherence theorem for monoidal categories 
(see \cite[Chapter~VII]{CWM}) it could be avoided. We write 
$\pi_{m,n}:X^mX^n\cong X^{m+n}$ for the canonical isomorphism built up
out of the associativity isomorphisms. (If \C really were strict this would
be the identity; otherwise, in order to make sense of tensor powers such
as $X^n$ some particular bracketing must be chosen.) 
We sometimes write $X$ for the
identity $1_X$ on an object $X$, and row vector notation $(f~g):A+B\to C$
for morphisms out of a coproduct. The composite of $f:X\to Y$ and $g:Y\to Z$
is written $g.f$.

\section{The approach of Dubuc}

The construction of a free monoid can be broken down into two parts. An
object $Y$ is said to be {\em pointed} if it is equipped with a map $y:I\to Y$;
we write \PtC for the category of {\em pointed objects} in \C.
Then the forgetful functor $U:\MonC\to\C$
is the composite of $V:\MonC\to\PtC$ which forgets the multiplication of a 
monoid but remembers the unit, and $W:\PtC\to\C$, which forgets the point. 
Since adjunctions compose, to find a left adjoint to $U=WV$, it will suffice
to find adjoints to $V$ and to $W$. But $W$ has a left adjoint sending 
$C\in\C$ to the coproduct injection $I\to I+C$, provided that coproducts with
$I$ exists, so in this case we are reduced to finding a left adjoint to $V$.
This reduction played a key role in \cite{Dubuc:free-monoids}, which 
contained a construction that will be important below (as well as various
transfinite variants which will not). We describe below one point of view
(not contained in \cite{Dubuc:free-monoids}) on this construction.

Thus we seek a left adjoint to $V:\MonC\to\PtC$.
In order to motivate the construction, we recall here
the connection between monoids and the simplicial category
\cite[Chapter~VII]{CWM}. We follow Mac~Lane in writing \DD for the category of 
finite ordinals and order-preserving maps: this is the ``algebraist's simplicial
category'', as opposed to the ``topologist's simplicial category'' which 
omits the empty ordinal, and reindexes the remaining objects. Now \DD is
monoidal with respect to ordinal sum, and ``classifies monoids in monoidal
categories'', in the sense that for any monoidal category
$(\C,\ot,I)$, the category \MonC of monoids in \C is equivalent to the 
category $\M(\DD,\C)$
of strong monoidal (=tensor-preserving) functors from \DD to \C, 
and monoidal natural
transformations. The strong monoidal functor corresponding to a monoid
$M$ in \C with multiplication $\mu:M^2\to M$ and unit $\eta:I\to M$ has image
$$\xymatrix{
I \ar[r]^{\eta} & M \ar@<3ex>[r]^{\eta M} \ar@<-3ex>[r]^{M\eta} & 
M^2 \ar[l]_{m} \ar@<6ex>[r]^{\eta M^2} \ar[r]^{M\eta M} \ar@<-6ex>[r]^{M^2\eta} &
M^3 \ar@<3ex>[l]_{\mu M} \ar@<-3ex>[l]_{M\mu} \ldots }$$
There is an analogous description of \PtC: let \DDm be the
(non-full) subcategory of \DD containing all the objects but only the 
injective order-preserving maps. This is still monoidal under ordinal sum,
and now \PtC is equivalent to 
the category $\M(\DDm,\C)$ of strong monoidal functors from $\DDm$ to \C 
and monoidal natural transformations. Corresponding to the pointed object
$(Y,y:I\to Y)$ we have
\begin{equation}\tag{$*$}
\xymatrix{
I \ar[r]^{y} & Y \ar@<2ex>[r]^{yY} \ar@<-2ex>[r]^{Yy} & 
Y^2 \ar@<3ex>[r]^{yY^2} \ar[r]^{YyY} \ar@<-3ex>[r]^{Y^2y} &
Y^3 \ldots 
 }
\end{equation}
If we identify \MonC with $\M(\DD,\C)$ and \PtC with $\M(\DDm,\C)$, then
the forgetful $V:\MonC\to\PtC$ is identified with the functor  
$\M(H,\C):\M(\DD,\C)\to\M(\DDm,\C)$ given
by composition with the inclusion $H:\DDm\to\DD$. If we were dealing with
ordinary functors rather than strong monoidal ones, in other words if we
sought a left adjoint to $\Cat(H,\C):\Cat(\DD,\C)\to\Cat(\DDm,\C)$, then
we could simply take the left Kan extension along $H$. In general this 
left Kan extension will not send strong monoidal functors to strong monoidal 
functors, but in special cases it does, and in fact provides the left
adjoint to $\M(H,\C)$. In such a case, if we form the 
strong monoidal functor $\DDm\to\C$ corresponding to a pointed object $(Y,y)$, 
and take its left Kan extension along $H$, the resulting strong monoidal 
functor from \DD to \C will correspond, via the equivalence 
$\M(\DD,\C)\simeq\MonC$, to a monoid in \C; and the underlying object of
this monoid is precisely the colimit of the diagram $(*)$ above, as a simple
calculation involving the coend formula for left Kan extensions shows. 
This, then, is the construction of Dubuc (in its simplest form where 
transfinite constructions are not required): if the
colimit of $(*)$ exists and is preserved by tensoring on either side, then
it has a monoid structure which is free on the pointed object $(Y,y:I\to Y)$.

\section{The  construction}

Colimits indexed by \DDm can be constructed iteratively using coequalizers
and colimits of chains, as we shall do below. Now many important functors 
do not preserve all coequalizers, but do preserve coequalizers of reflexive 
pairs (pairs which have a common section). There is also a general way to 
replace a pair 
$f,g:A \rightrightarrows B$ by a reflexive pair with the same coequalizer: 
replace $A$ by $A+B$, and then use the identity map on $B$, as in 
$$\xymatrix{
A+B \ar@<1ex>[r]^-{(f~B)} \ar@<-1ex>[r]_-{(g~B)} & B. }$$
The construction given here
amounts to  an analogous adaptation of the construction of Dubuc described
in the previous section. 
Of course when both constructions work, they agree; in particular this will
be the case if coproducts are preserved by tensoring on either side.

Suppose then that $(Y,y:I\to Y)$ is a pointed object, and form the 
corresponding diagram $(*)$. The colimit of $(*)$ can be constructed as
follows. For each $n$ and each $k=0,\ldots n-2$, form the coequalizer
$$\xymatrix{
Y^{n-1} \ar@<1ex>[rr]^-{Y^kyY^{n-k-1}} 
\ar@<-1ex>[rr]_-{Y^{k+1}yY^{n-k-2}} &&
Y^n \ar[r]^-{r_{n,k}} & Y^n_k. }$$
Now form the cointersection $r_n:Y^n\to Y_n$ of the $r_{n,k}:Y^n\to Y^n_k$
as $k$ runs from $0$ to $n-2$. A straightforward calculation shows that for
all $j$ and $k$, the composites
$$\xymatrix{
Y^{n-2} \ar@<1ex>[rr]^{Y^k y Y^{n-k-2}} \ar@<-1ex>[rr]_{Y^{k+1}y Y^{n-k-3}} &&
Y^{n-1} \ar[rr]^{Y^j y Y^{n-j-1}} && Y^n \ar[r]^{r_n} & Y_n }$$
agree, while the composite $r_n.Y^j y Y^{n-j-1}$ is independent of $j$, 
and so by the universal property of $r_{n-1}:Y^{n-1}\to Y_{n-1}$ there is 
a unique map $h_n:Y_{n-1}\to Y_n$ such that the square
$$\xymatrix{
Y^{n-1} \ar[rr]^{Y^jyY^{n-j-1}} \ar[d]_{r_{n-1}} && Y^n \ar[d]^{r_n} \\
Y_{n-1} \ar[rr]_{h_n} && Y_n }$$
commutes
for all $j$. The colimit of $(*)$, corresponding to the construction of Dubuc,
is the colimit of the chain consisting of the $Y_n$ with connecting maps $h_n$.

We modify this at the first step only, replacing coequalizers by reflexive
coequalizers, as follows. The original coequalizers can be written as
$$\xymatrix{
Y^{n-1}=Y^kYY^{n-k-2} \ar@<1ex>[rr]^-{Y^kyYY^{n-k-2}} 
\ar@<-1ex>[rr]_-{Y^kYyY^{n-k-2}} && 
Y^n \ar[r]^-{r_{n,k}} & Y^n_k }$$
and our new reflexive coequalizers are
$$\xymatrix {
Y^k(Y+Y^2)Y^{n-k-2} \ar@<3ex>[rr]^-{f_{n,k}}
\ar@<-3ex>[rr]^-{g_{n,k}} && 
Y^n \ar[ll]_-{d_{n,k}} \ar[r]^{q_{n,k}} & Z^n_k }$$
where $f_{n,k}=Y^k(yY~Y^2)Y^{n-k-2}$, 
$g_{n,k}=Y^k(Yy~Y^2)Y^{n-k-2}$, and $d_{n,k}=Y^kbY^{n-k-2}$, with $b$ the 
coproduct injection $Y^2\to Y+Y^2$. 
Then, as before, we shall form the cointersection $q_n:Y^n\to Z_n$ of the
$q_{n,k}$, the induced maps $j_n:Z_{n-1}\to Z_n$ (see below) satisfying 
$j_n. q_{n-1}=q_n. (Y^kyY^{n-k-1})$, and the colimit $Z$ of the chain 
consisting of the $Z_n$ and the $j_n$; we write 
$z_n:Z_n\to Z$ for the legs of the colimit cocone. 

It is, however, worth taking a little more time to justify the existence
of the $j_n$. We must show that for all $j$ and $k$, the composites
$$\xymatrix{
Y^k(Y+Y^2)Y^{n-k-3} \ar@<1ex>[rr]^-{Y^k(yY~Y^2)Y^{n-k-3}} 
\ar@<-1ex>[rr]_-{Y^k(Yy~Y^2)Y^{n-k-3}} && Y^{n-1} 
\ar[rr]^-{Y^j y Y^{n-j-1}} && Y^n \ar[r]^{q_n} & Z_n}$$
are equal. This is completely straightforward if either $j\le k$ or 
$j\ge k+2$, but the case $j=k+1$ is a bit more complicated; it can be 
broken down as in the following diagram:

$$
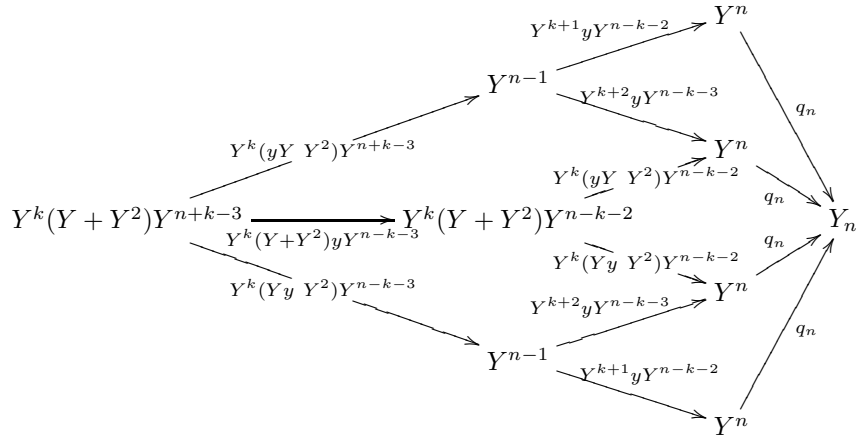
\xymatrix @R1pc {
&&& Y^n \ar[dddr]^{q_n} \\
&& Y^{n-1} \ar[ur]^{Y^{k+1}y Y^{n-k-2}} \ar[dr]^{Y^{k+2}y Y^{n-k-3}} \\
&&& Y^n \ar[dr]_{q_n} \\
Y^k(Y+Y^2)Y^{n+k-3} \ar[uurr]|{Y^k(yY~Y^2)Y^{n+k-3}} 
\ar[rr]_{Y^k(Y+Y^2)yY^{n-k-3}} \ar[ddrr]|{Y^k(Yy~Y^2)Y^{n-k-3}} &&
Y^k(Y+Y^2)Y^{n-k-2} \ar[ur]|(0.6){Y^k(yY~Y^2)Y^{n-k-2}} 
\ar[dr]|(0.6){Y^k(Yy~Y^2)Y^{n-k-2}} && Y_n \\
&&& Y^n \ar[ur]^{q_n} \\
&& Y^{n-1} \ar[ur]^{Y^{k+2}y Y^{n-k-3}} \ar[dr]^{Y^{k+1}y Y^{n-k-2}} \\
&&& Y^n \ar[uuur]_{q_n} 
}$$
in which the individual regions are easily seen to commute.

In the following section we prove that $Z$ can be made into a monoid
which is free on $(Y,y)$. We record the general result as:

\begin{theorem}
Let \C be a monoidal category with finite limits and countable colimits,
and the functors $-\ot C$ and $C\ot -$ preserve reflexive coequalizers and
colimits of countable chains. This includes in particular the case where \C 
has the stated limits and colimits and $C\ot-$ and $-\ot C$ preserve sifted
colimits. Then the free monoid on a pointed object
$(Y,y)$ exists, and its underlying object $Z$ can be calculated as above. 
The free monoid on an object $X$ is found by taking $Y$ to be $I+X$ and $y$ 
to be the coproduct injection. 
\end{theorem}

\section{The proof}

Suppose that \C satisfies the conditions of the theorem.
Our construction involved three
types of colimit: reflexive coequalizers, finite cointersections of regular
epimorphisms, and colimits of chains. By assumption, the first and third of
these are preserved by tensoring; we shall see that the second is also 
preserved. We defer for the moment the proof, merely noting that the 
case of binary cointersections suffices, and recording:

\begin{lemma}\label{lemma:cointersection}
If $q:B\to C$ and $q':B\to C'$ are regular epimorphisms, then their
cointersection (pushout)   
$$\xymatrix{
B \ar[r]^{q} \ar[d]_{q'} & C \ar[d]^{r} \\
C' \ar[r]_{r'} & D }$$
is preserved by tensoring on either side. 
\end{lemma}

We also need the following form of the ``3-by-3 lemma''
\cite{ptj:topos-theory} for reflexive 
coequalizers, whose proof is once again deferred.

\begin{lemma}[3-by-3 lemma]\label{lemma:3by3}
If 
$$\xymatrix{
A_1 \ar@<2ex>[r]^{h_1}\ar@<-2ex>[r]_{h_2} & A_2 \ar[l] \ar[r]^{h} & A_3 \\
B_1 \ar@<2ex>[r]^{k_1}\ar@<-2ex>[r]_{k_2} & B_2 \ar[l] \ar[r]^{k} & B_3 
}$$
are reflexive coequalizers, preserved by tensoring on either side, then
$$\xymatrix{
A_1\ot B_1 \ar@<2ex>[r]^{h_1\ot k_1} \ar@<-2ex>[r]_{h_2\ot k_2} & 
A_2\ot B_2 \ar[l] \ar[r]^{h\ot k} & A_3\ot B_3 }$$
is also a reflexive coequalizer, and $A_3\ot B_3$ is the cointersection 
of $A_3\ot B_2$ and $A_2\ot B_3$ (as quotients of $A_2\ot B_2$).
This shows in particular that regular epimorphisms are closed under
tensoring. 
\end{lemma}

We need to construct a multiplication $\mu:ZZ\to Z$. The idea will be
first to construct $\mu_{m,n}:Z_mZ_n\to Z_{m+n}$, then show that they
pass to the colimit to give the desired $\mu$. 

\subsection{Construction of $\mu_{m,n}$}

Consider first the diagram
$$\xymatrix{
Y^mZ^n_l \ar[r]^{\cong} & Z^{m+n}_{m+l}  \\
Y^mY^n \ar[u]^{Y^m q^n_l} \ar[d]_{q^m_k Y^n} \ar[r]^{\pi_{m,n}} & 
Y^{m+n} \ar[u]^{q^{m+n}_{m+l}} \ar[d]_{q^{m+n}_k} \\
Z^m_k Y^n \ar[r]^{\cong} & Z^{m+n}_k 
}$$
which we shall build out of the canonical isomorphism 
$\pi_{m,n}:Y^mY^n\cong Y^{m+n}$.
Since $q^m_k:Y^m\to Z^m_k$ was constructed as the coequalizer of maps 
$f_{m,k}$ and $g_{m,k}$, thus $q^m_kY^n:Y^{m+n}=Y^mY^n\to Z^m_kY^n$ can
be constructed as the coequalizer of $f^m_kY^n$ and $g^m_kY^n$; that is,
of $f^{m+n}_k$ and $g^{m+n}_k$; thus we get the induced isomorphism 
$Z^m_kY^n\cong Z^{m+n}_k$ at the bottom of the diagram. Similarly the
coequalizer defining $Z^n_l$ is preserved by tensoring on the left by $Y^m$
and so we get the induced isomorphism $Y^mZ^n_l\cong Z^{m+n}_{m+l}$ at the
top. Thus $Z_mY^n$ is the cointersection of all the $Z^{m+n}_p$ with 
$0\le p\le m-2$, and $Y^mZ_n$ is the cointersection of all the $Z^{m+n}_p$
with $m\le p\le m+n-2$. By the 3-by-3 lemma,
$$\xymatrix{
Y^mY^n \ar[r]^{Y^mq_n} \ar[d]_{q_m Y^n} & Y^m Z_n \ar[d]^{q_m Z_n} \\
Z_m Y^n \ar[r]_{Z_m q_n} & Z_mZ_n }$$
is a cointersection, and so $Z_mZ_n$ is the cointersection of all the 
$Z^{m+n}_p$ with $0\le p\le m-2$ or $m\le p\le m+n-2$. On the other hand
$Z_{m+n}$ is the cointersection of all the $Z^{m+n}_p$ with $0\le p\le m+n-2$,
and so there is a canonical quotient map $\mu_{m,n}:Z_mZ_n\to Z_{m+n}$ fitting
into the commutative diagram
$$\xymatrix{
Y^mY^n \ar[r]^{\pi_{m,n}} \ar[d]_{q_m q_n} & Y^{m+n} \ar[d]^{q_{m+n}} \\
Z_m Z_n \ar[r]_{\mu_{m,n}} & Z_{m+n} }$$

\subsection{Construction of $\mu$}

Since tensoring preserves colimits of chains, we have
$$ZZ=Z\ot Z=(\colim_m Z_m)\ot(\colim_n Z_n)\cong\colim_m\colim_n(Z_m\ot Z_n)$$
so there will be a unique map $\mu:ZZ\to Z$  making
$$\xymatrix{
Z_mZ_n \ar[r]^{\mu_{m,n}} \ar[d] & Z_{m+n} \ar[d] \\
ZZ \ar[r]_{\mu} & Z }$$
commute provided that the $\mu_{m,n}$ are compatible with the maps 
$j_n:Z_n\to Z_{n+1}$ (and $j_m$); in other words that the maps $\mu_{m,n}$ 
are natural in $m$ and $n$. We explain the naturality in $n$; the case of
$m$ is similar. In the first diagram below, the left square commutes by 
definition of $\mu_{m,n}$, and the right square commutes by definition 
of $j_{m+n+1}$. In the second diagram, commutativity of the left square
follows from the definition of $j_{n+1}$, while the right square commutes
by definition of $\mu_{m,n+1}$.
$$
\xymatrix{
Y^mY^n \ar[r]^{\pi_{m,n}} \ar[d]_{q_m q_n} & Y^{m+n} \ar[d]_{q_{m+n}}
\ar[r]^{Y^{m+n}y} & Y^{m+n+1} \ar[d]_{q_{m+n+1}} \\
Z_m Z_n \ar[r]_{\mu_{m,n}} & Z_{m+n} \ar[r]_{j_{m+n+1}} & Z_{m+n+1}
}\xymatrix{
Y^mY^n \ar[r]^{Y^mY^ny} \ar[d]_{q_m q_n} & Y^m Y^{n+1} \ar[d]_{q_m q_{n+1}} 
\ar[r]^{\pi_{m,n+1}} & Y^{m+n+1} \ar[d]_{q_{m+n+1}} \\
Z_m Z_n \ar[r]_{Z_m j_{n+1}} & Z_m Z_{n+1} \ar[r]_{\mu_{m,n+1}} & Z_{m+n+1} }
$$
Now the composites across the top of the two diagrams agree, by naturality
of associativity, and the (common) left vertical $q_m q_n$ is a regular 
epimorphism, so that the composites across the bottom agree. This gives
the desired naturality in $n$, and so we obtain the required map 
$\mu:ZZ\to Z$.

\subsection{Verification of associative and unit laws}

The associative law $\mu. \mu Z=\mu. Z\mu$ will hold provided
that
$$\xymatrix{
Z_mZ_nZ_p \ar[r]^{\mu_{m,n} Z_p} \ar[d]_{Z_m \mu_{n,p}} & 
Z_{m+n}Z_p \ar[d]^{\mu_{m+n,p}} \\
Z_mZ_{n+p} \ar[r]_{\mu_{m,n+p}} & Z_{m+n+p} }$$
commutes for all $m$, $n$, and $p$. Now the two paths around this square
will agree provided that they agree when composed with the regular
epimorphism $q_mq_nq_p:Y^mY^nY^p\to Z_mZ_nZ_p$, and this in turn follows
from the evident commutativity of 
$$\xymatrix{
Y^mY^nY^p \ar[r]^{\pi_{m,n} Y^p} \ar[d]_{Y^m \pi_{n,p}} & 
Y^{m+n}Y^p \ar[d]^{\pi_{m+n,p}} \\
Y^mY^{n+p} \ar[r]_{\pi_{m,n+p}} & Y^{m+n+p}. }$$

The unit is given by the composite
$$\xymatrix{
I \ar[r]^-{y} & Y = Z_1 \ar[r]^-{z_1} & Z}$$
where $z_1$ is the relevant leg of the colimit cocone; the
verification of the unit law is similar to but easier than the verification
of associativity.

\subsection{Universal property}

The unit of the adjunction will be the map
$$\xymatrix{
Y=Z_1 \ar[r]^-{z_1} & Z}$$
of pointed objects; we must show that this has the appropriate universal 
property. In other words, for every monoid $M=(M,\mu,\eta)$ and every 
morphism $f:(Y,y)\to(M,\eta)$ of pointed objects, we must show that there
is a unique monoid morphism $g:(Z,\mu,\eta)\to(M,\mu,\eta)$ with $gz_1=f$. 

For each $n$, we have the composite $f_n$ as in 
$$\xymatrix{
Y^n \ar[r]^{f^n} & M^n \ar[r]^{\mu_{(n)}} & M }$$
where $\mu_{(n)}$ is the $n$-ary multiplication operation for the monoid $M$.
We must show that these maps $f_n=\mu_{(n)}. f^n$ pass to the quotient to 
give $g_n:Z_n\to M$. We check only that the composites in
$$\xymatrix{
Y+Y^2 \ar@<1ex>[rr]^{(yY~Y^2)} \ar@<-1ex>[rr]_{(Yy~Y^2)} && 
Y^2 \ar[r]^{f^2} & M^2 \ar[r]^{\mu} & M }$$
are equal; the other cases all follow by functoriality of $\ot$. Now the
two displayed composites are maps out of a coproduct, so will agree if 
their components do; for the components on $Y^2$ this is trivial, and for
the components on $Y$ we have
$$\mu. f^2. yY=\mu. \eta M. f=f=\mu. M\eta. f=
\mu. f^2. Yy.$$
Thus the maps $\mu_{(n)}. f^n$ induce maps $g_n:Z_n\to M$, which clearly
pass to the colimit to give $g:Z\to M$. We must show that this is a monoid
map, and is the unique such which extends $f$. 

Now $g$ preserves the unit by construction, and will preserve the multiplication
provided that 
$$\xymatrix{
Z_mZ_n \ar[r]^{g_m g_n} \ar[d]_{\mu_{m,n}} & MM \ar[d]^{\mu} \\
Z_{m+n} \ar[r]_{g_{m+n}} & M }$$
commutes. But $Z_mZ_n$ is a quotient of $Y^mY^n$, so this in turn restricts
to commutativity of 
$$\xymatrix{
Y^mY^n \ar[r]^{f_m f_n} \ar[d]_{\pi_{m,n}} & MM \ar[d]^{\mu} \\
Y^{m+n} \ar[r]_{f_{m+n}} & M}$$
which holds by construction of the $f_n$ and associativity of $\mu$. 

This proves that $g$ is a monoid map; it remains to show the uniqueness.
Suppose then that $h:Z\to M$ is a monoid map, with $h. z_1=f$. In order to 
show that $h=g$, it will suffice to show that $h. z_n=g_n$ for all $n$. 
This in turn will hold if $h. z_n. q_n=f_n$ for all $n$. Thus we 
must show that the exterior of the diagram
$$\xymatrix{
Y^n \ar[r]^{z^n_1} \ar[d]_{q_n} & Z^n \ar[r]^{h^n} \ar[d]^{\mu_{(n)}} & 
M^n \ar[d]^{\mu_{(n)}} \\
Z_n \ar[r]_{z_n} & Z \ar[r]_{h} & M }$$
commutes. The right square commutes because $h$ is a monoid homomorphism,
so it suffices to show that the left square commutes, and this follows from
the definition of $\mu:Z^2\to Z$ by a straightforward induction. 

\subsection{Proof of lemmas}

Consider a diagram
$$\xymatrix{
A_{11} \ar@<1ex>[r]^{f_1} \ar@<-1ex>[r]_{f_2} 
\ar@<-1ex>[d]_{f'_2} \ar@<1ex>[d]^{f'_1}
& A_{12} \ar@<-1ex>[d]_{g'_2} \ar@<1ex>[d]^{g'_1} \\
A_{21} \ar@<1ex>[r]^{g_1} \ar@<-1ex>[r]_{g_2} 
& A_{22} }$$
in which $g_i. f'_j=g'_j. f_i$ for $i,j\in\{1,2\}$, and suppose
also that there exist $s:A_{12}\to A_{11}$ and $s':A_{21}\to A_{11}$ with
$f_1. s=f_2. s=1$ and $f'_1. s'=f'_2. s'=1$. 

Then a map $x:A_{22}\to B$ satisfies $x. g'_1. f_1=x. g'_2. f_2$
if and only if it satisfies $x. g'_1=x. g'_2$ and 
$x. g_1=x. g_2$. For if the former equation holds then we have
$$x. g'_1=x. g'_1. f_1. s=x. g'_2. f_2. s=x. g'_2$$
$$
x. g_1=x. g_1. f'_1. s'=x. g'_1. f_1. s'=
x. g'_2. f_2. s'=x. g_2. f'_2. s'=x. g_2$$
\noindent while if the latter two equations hold then
$$x. g'_1. f_1=x. g'_2. f_1=x. g_1. f'_2=
x. g_2. f'_2=x. g'_2. f_2.$$
As a result we have:

\begin{proposition}
In the situation above, the coequalizer of  $g'_1. f_1$ and $g'_2. f_2$
is the cointersection of the coequalizer of $g_1$ and $g_2$ and the 
coequalizer of $g'_1$ and $g'_2$. 
\end{proposition}

To prove the 3-by-3 lemma (Lemma~\ref{lemma:3by3}), apply this in the 
case of 
$$\xymatrix{
A_1\ot B_1 \ar@<1ex>[r]^{A_1\ot k_1} \ar@<-1ex>[r]_{A_1\ot k_2} 
\ar@<-1ex>[d]_{h_2\ot B_1} \ar@<1ex>[d]^{h_1\ot B_1}
& A_1\ot B_2 \ar@<-1ex>[d]_{h_2\ot B_2} \ar@<1ex>[d]^{h_1\ot B_2} \\
A_2\ot B_1 \ar@<1ex>[r]^{A_2\ot k_1} \ar@<-1ex>[r]_{A_2\ot k_2} 
& A_2\ot B_2 }$$
noting that the $A_1\ot k_i$ have a common section $A_1\ot t$, and the
$h_1\ot B_1$ have a common section $s\ot B_1$. 

To prove Lemma~\ref{lemma:cointersection}, let $q:B\to C$ and $q':B\to C'$
be the coequalizers of the reflexive pairs
$$\xymatrix{
A \ar@<2ex>[r]^{h_1} \ar@<-2ex>[r]_{h_2} & B \ar[l]_{s} &
A' \ar@<2ex>[r]^{h'_1} \ar@<-2ex>[r]_{h'_2} & B \ar[l]_{s'} &
}$$
and   form the universal object $P$ with morphisms
$$\xymatrix{
P \ar@<1ex>[r]^{k_1} \ar@<-1ex>[r]_{k_2} 
\ar@<-1ex>[d]_{k'_2} \ar@<1ex>[d]^{k'_1}
& A' \ar@<-1ex>[d]_{h'_2} \ar@<1ex>[d]^{h'_1} \\
A \ar@<1ex>[r]^{h_1} \ar@<-1ex>[r]_{h_2} 
& B }$$
satisfying equations as above. In terms of elements, this would be formed
as $\{x_1,x_2\in A, x'_1,x'_2\in A'\mid  h_i(x'_j)=h'_j(x_i) \}.$ It is 
straightforward to show that the relevant pairs are reflexive; for example
$x'\mapsto \bigl(sh'_1(x'),sh'_2(x'),x',x'\bigr)$ provides a common 
section to $k_1$ and $k_2$. 
The proposition then reduces the cointersection of $q$ and $q'$ to 
the reflexive coequalizer of $h'_1.k_1$ and $h'_2.k_2$, which 
by assumption is preserved by tensoring.

\section{Examples}\label{sect:examples}

If \C is any variety, equipped with the cartesian product $\times$
as tensor product, then products, reflexive coequalizers, and colimits
of chains are all computed as in \Set, and since the product in \Set
with a fixed object is cocontinuous, it follows that tensoring in \C with
a fixed object preserves the relevant colimits.

For a fixed set $A$, the category $\Span(A,A)$ of spans from $A$ to $A$ is
the category of all sets over $A\times A$. This is monoidal via pullback
$$
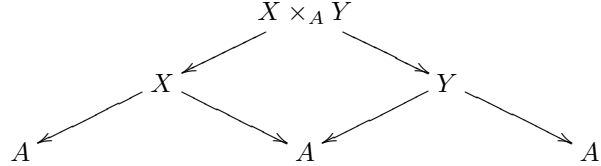
\xymatrix @!R @!C @R1pc @C1pc {
&& X\times_A Y \ar[dl] \ar[dr] \\
& X \ar[dl] \ar[dr] && Y \ar[dl] \ar[dr] \\
A && A && A }$$
and a monoid in the resulting monoidal category is precisely a category
with object-set $A$.
Tensoring on either side is cocontinuous (and in fact has an adjoint) 
because pullbacks in \Set are cocontinuous. But now we can move from \Set
to a category \E with finite limits in which pullback may not be cocontinuous,
but does preserve reflexive coequalizers and colimits of chains: this is
the case, for example, in any variety. If we consider an object
$A\in\E$, and the category $\Span(\E)(A,A)$ of internal spans in \E from 
$A$ to $A$ this is once again monoidal, and once again a monoid is a
category in \E with $A$ as its object of objects. But this time tensoring
on either side preserves reflexive coequalizers and colimits of chains,
but not arbitrary colimits. Thus our construction gives free internal categories
in \E.

For a more structured example, one could consider not $\Span(\E)$ but
$\Prof(\E)$, the bicategory of internal categories in \E and profunctors 
between them. Fixing an internal category $A$, we get a monoidal category
$\Prof(\E)(A,A)$, and once again the conditions for our construction 
will be satisfied. Taking $\E=\Mon$, the category of monoids, and $A$ to 
be (a suitable strict version of) the monoidal category \PP of finite 
sets and bijections, we get a monoidal category $\Prof(\Mon)(\PP,\PP)$ 
in which monoids are precisely PROPs (see \cite{prop}), and so a different
notion of free PROP to that given in \cite{Vallette:free-monoids}. 

Finally, for a slightly childish example, take the category \C to be the 
category \Grp of groups and group homomorphisms, with the cartesian monoidal 
structure (with the product as tensor product). For a group $G$ the functors
$G\times-:\Grp\to\Grp$ and $-\times G:\Grp\to\Grp$ do not of course preserves all
colimits, but they do preserve reflexive coequalizers and colimits of chains,
as would be the case with any variety in place of \Grp. Now by the ``Eckmann-Hilton argument'', a monoid in \Grp is precisely an abelian group. So our 
construction reduces to the abelianization of a group. 

\section{Algebraically free monoids}
\label{sect:algebraically-free}

The free monoid construction we have been looking at involves a
pointed object $(Y,y)$ a monoid $(Z,\mu,\eta)$, and a morphism of
pointed objects $k:(Y,y)\to(Z,\eta)$. But there is another possible
universal property that such data might satisfy. Write $\C^{(Z,\mu,\eta)}$
for the category of objects of \C equipped with an action of the monoid
$(Z,\eta,\mu)$, and write $\C^{(Y,y)}$ for the category of objects of \C
equipped with an action of the pointed object $(Y,y)$: in other words,
a morphism $\alpha:YA\to A$ satisfying $\alpha.yA=1$. There is a functor
$k^*:\C^{(Z,\eta,\mu)}\to\C^{(Y,y)}$ sending $A$ equipped with $Z$-action 
$\beta:ZA\to A$ to $A$ equipped with $Y$-action 
$$\xymatrix{YA \ar[r]^{kA} & ZA \ar[r]^{\beta} & A.}$$
When this functor $k^*$ is an isomorphism of categories we say that 
$k$ exhibits $(Z,\mu,\eta)$ as the {\em algebraically free} monoid on 
$(Y,y)$ \cite{Kelly-transfinite}. The algebraically free monoid, if it
exists, is free, but it is possible for a free monoid to exist without
being algebraically free; nonetheless under fairly general conditions the 
algebraically free monoid exists (and is free); see 
\cite[Section~23]{Kelly-transfinite}. 
This means that the free monoid on $(Y,y)$ can be found by calculating
free $(Y,y)$-actions, an idea that implicitly goes back to  
\cite{Barr:free-triples}. All cases where the free monoid is computed in
\cite{Kelly-transfinite} are done in this way. 

Under the hypotheses of our theorem, a
necessary and sufficient condition for the free monoid on
$(Y,y)$ to be algebraically free is that for
any action  $\alpha:YA\to A$ of $(Y,y)$ on $A$, the composites
$$\xymatrix{
(Y+Y^2)A \ar@<1ex>[rr]^-{(yY~Y^2)A} \ar@<-1ex>[rr]_-{(Yy~Y^2)A} &&
Y^2A \ar[r]^{Y\alpha} & YA \ar[r]^{\alpha} & A }$$
agree. 

In general there seems to be no reason why this should always be true, 
although we do not have a specific example where it fails. We therefore
conjecture that the hypotheses of our theorem are not sufficient to 
guarantee that the free monoid is algebraically free.

\bibliographystyle{plain}

\end{document}